\documentclass[12pt]{article}

\usepackage{hyperref}
\usepackage[latin1]{inputenc}
\usepackage{makeidx}
\usepackage[french]{babel}
\usepackage{amsfonts}
\usepackage{amssymb}
\usepackage{amsthm}
\usepackage{newlfont}
\usepackage{showidx}
\usepackage[dvips]{graphicx}

\newtheorem{theo}{Théorème}

\newtheorem{cor}{Corollaire}
\newtheorem{rem}{Remarque}

\newtheorem{lem}{Lemme}

\def\ad{\mathrm{ad}}
\def\Res{\mathrm{Res}}
\def\Rad{\mathrm{Rad}}

\newcommand{\C}{\mathbb{C}}
\newcommand{\N}{\mathbb{N}}
\def\ho{ho\-mo\-mor\-phis\-me}

\def\A{\alpha}
\def\D{\Delta}

\def\p{\mathfrak{p}}
\def\k{\mathfrak{k}}
\def\m{\mathfrak{m}}
\def\aa{\mathfrak{a}}
\def\g{\mathfrak{g}}
\def\h{{\mathfrak h}}
\def\s{{\mathfrak s}}
\def\sf{\s_{f}}
\def\gf{\g(f)}

\def\go{\g_{o}}

\begin{document}

\bibliographystyle{alpha}
\author{Charles Torossian\footnote{UMR 8553 du CNRS, DMA-ENS, École Normale Supérieure, 45 rue d'Ulm 75230 Paris
cedex 05, Charles.Torossian@ens.fr}}
\title{Deux résultats autour du théorème de restriction de Chevalley}\maketitle

{\noindent \bf{Résumé: }}Dans ce court article, nous exposons deux
résultats. Le premier donne une preuve nouvelle de la surjectivité
dans le théorème de restriction de Chevalley; elle se base sur
l'utilisation des opérateurs de Dunkl.  Le second résultat décrit
précisément l'image par l'homomorphisme de restriction dans le cas
des  algèbres de Takiff.
On étend ces résultats pour les espaces symétriques de Takiff.\\

{\noindent \bf{Abstract:}}  We use Dunkl's operators to give an
elementary proof of the surjectivity in the Chevalley's restriction
theorem. In the second part of this article we describe the image of
the invariants by the restriction map  in the case of Takiff
algebras. We extend this result to Takiff symmetric spaces.
\\

\section{Espaces symétriques et opérateurs de Dunkl}
\subsection{Opérateurs de Dunkl} Soit $G$ un groupe de Lie
semi-simple r\'eel \`a centre fini et $K$ un sous-groupe compact
maximal. On note $\g$ l'alg\`ebre de Lie de $G$, $\k$ celle de $K$
et
 $\g=\k \oplus \p$ la d\'ecomposition de Cartan associ\'ee. Soit
$\aa$ un sous espace de Cartan dans $\p$, $\Delta$ le syst\`eme de
racines associ\'e, $\Sigma$ le syst\`eme des racines indivisibles de
$\Delta$ et $W$ le groupe de Weyl.
 On note
$S[\p]$ et $S[\aa]$ les algèbres symétriques de $\p$ et $\aa$. On
note $S[\p]^K$ les éléments $K-$invariants de $S[\p]$ et $S[\aa]^W$
les éléments $W-$invariants de $S[\aa]$. On notera par $S(\p)$
(resp. $S(\aa)$) le corps des fractions de $S[\p]$ (resp. $S(\aa)$)
et par $S(\p)^K$ (resp. $S(\aa)^W$) les fractions
$K$-invariantes (resp. $W$-invariantes).\\


La décomposition, via la forme de Killing, $\p=\aa\oplus \aa^\perp$
permet de définir l'homomorphisme de restriction de $S[\p]^K$ dans
$S[\aa]^W$. \\

Pour $\alpha \in\Sigma$ on note $r_{\alpha}$ la r\'eflexion par
rapport \`a la coracine $H_{\alpha}$ c'est à dire que l'on a
$$r_{\alpha}(H)=H-\frac{2 \alpha(H) }{\alpha(H_{\alpha})}H_{\alpha}.$$ On
note aussi $k_{\alpha}= \frac{1}{2}(m_{\alpha} + m_{2 \alpha})$ avec
$m_{\alpha}$  la multiplicit\'e de la racine $\alpha$.

\noindent Pour $\xi \in \aa$ on définit  $T_{\xi}$ l'op\'erateur de
Dunkl \cite{dunkl1} sur $\aa$ par :
 $$T_{\xi}= \partial_{\xi} + \frac{1}{2}\sum_{\alpha \in \Sigma} k_{\alpha} \alpha ( \xi)
\frac{1-r_{\alpha}}{\alpha (.)}.$$Alors $T_{\xi}$ est bien d\'efini
sur les fonctions polynomiales sur $\aa$; il est de degr\'e~$-1$.
Ces op\'erateurs ont des propri\'et\'es remarquables que nous
r\'esumons (cf. \cite{dunkl1},~\cite{heck}).
\begin{enumerate}
\item Ils commutent deux \`a deux, ce qui permet de d\'efinir $p(T)$
 pour $p$ dans l'alg\`ebre
symétrique \ $S[\aa]$. Ce point est assez difficile et est dû à
Dunkl (d'autres démonstrations existent dues à Opdam et Heckman).
\item Pour $p \in S[\aa]^{W}$, l'op\'erateur $p(T)$ est un op\'erateur diff\'erentiel
dans l'espace des fonctions $W$-invariantes. Ce point est clairement
facile.
\item Pour $P$ un \'el\'ement de $S[\p]^{K}$ on note $\partial(P)$ l'op\'erateur
diff\'erentiel $K$-invariant sur $\p$ associ\'e, alors la partie
radiale de $\partial(P)$ est \'egale \`a l'action de $p(T)$ sur les
éléments $W$-invariants, o\`u $p$ est est la restriction de $P$ \`a
$\aa$. Ce point se démontre en vérifiant la formule pour le
laplacien $L$ sur $\p$ par un calcul direct, puis en utilisant la
formule suivante pour $P$ homog\`ene de degr\'e $n$
$$\frac{1}{2^{n} n!} \ad^{n}(L)(P^{\sharp})=\partial(P)$$
o\`u $P^{\sharp}$ est la fonction polynomiale sur $\p$ correspondant
\`a $P$ (via la forme de Killing). On en d\'eduit que l'on a
$$\Rad\partial(P)=\frac{1}{2^{n} n!}\ad^{n}(\ell(T))(p^{\sharp})$$
sur les fonctions $W-$ invariantes sur $\aa$ (on a not\'e
$\ell=\Res(L)$). Puis pour les m\^emes raisons on a
$$\frac{1}{2^{n} n!}\ad^{n}(\ell(T))(p^{\sharp})=p(T).$$

\end{enumerate}

\subsection{Surjectivité de la restriction de Chevalley}
\begin{theo}[Chevalley] L'application de restriction de $S[\p]^K$ sur $S[\aa]^W$ est un
isomorphisme. \end{theo}

La partie difficile du théorème est dans la surjectivité. En effet
il est facile de voir que la trace des $K$-orbites génériques sur
$\aa$ forment des $W$-orbites. Il en résulte d'une part
l'injectivité de l'application de restriction et d'autre part
l'isomorphisme par l'homomorphisme de restriction au niveau des
corps de fractions : $S(\p)^K \approx S(\aa)^W$ (cf. paragraphe
\ref{res} pour des considérations similaires
plus générales).\\

Voici un argument élémentaire utilisant les opérateurs de Dunkl.\\

Notons $R$ l'image de l'algèbre $S[\p]^K$ par l'application de
restriction de Chevalley. C'est une sous-algèbre de $S[\aa]^W$.\\

Il existe un produit scalaire (positif, non dégénéré) sur $S[\aa]$
donné par la formule $$\langle p, q \rangle= p(T)(q)(0).$$ Cette
formule définit encore un produit scalaire  pour toutes les
multiplicités positives (cf. \cite{jeu}). Les opérateurs de Dunkl
sont équivariants par rapport à l'action de $W$ et la restriction de
ce produit à $S[\aa]^W$ est encore non dégénérée et positive.

\begin{lem}
La restriction réalise une isométrie entre l'espace $S[\p]^K$ muni
du produit de Fischer et l'espace $R$ muni du produit scalaire de
Dunkl.
\end{lem}

\noindent \textit{Preuve :} Pour $P, Q$ dans $S[\p]^{K}$ il suffit
de remarquer que l'on a $$\langle p, q \rangle= P(\partial)(Q)(0),$$
avec $p,q$ les restrictions de $P,Q$. Cela résulte de la définition
de la partie radiale de l'opérateur $P(\partial)$ et de la
propriété~3 du paragraphe précédent. Comme le produit de Fischer
est non dégénéré et positif sur $S[\p]^K$ on en déduit le résultat.\\

Notons $R^{\perp}$ l'orthogonal de $R$ dans $S[\aa]^W$  pour le
produit scalaire de Dunkl. On a donc  $$R\oplus R^\perp= S[\aa]^W.$$
Il résulte de la définition de ce produit que  l'adjoint de $T_\xi$
($\xi \in \aa$) est l'opérateur de multiplication par $\xi$.  Par
ailleurs $R$ est stable pour l'action de $p(T)$ pour $p\in R$, car
$S[\p]^K$ est stable par l'action des opérateurs $\partial(P)$, on
en déduit que $R^\perp$
est stable pour la multiplication par $p\in R$: c'est le point clé de notre argument.\\

Tout élément de $S[\aa]^W$ peut être remonté en une fraction $K$
invariante sur $\p$. C'est le côté facile de l'isomorphisme de
Chevalley. On sait même que c'est une fraction bien définie près des
éléments semi-simples réguliers. Ces éléments sont définis par
l'inégalité $\pi\neq 0$ avec $\pi$ le polynôme $K-$invariant
coefficient du monôme de plus bas degré en $T$ obtenu dans le
développement de $\det_\p(\mathrm{ad}^2(X) -T)$.\\

Terminons notre preuve du théorème de Chevalley.\\

Soit $p\in R^\perp$ et notons $F$ la fraction $K-$invariante
correspondante. Alors pour tout $r \in S[\p]^K$ non nul, $r F$ est
une fraction qui n'est pas un polynôme, car $R^\perp$ est stable par
multiplication par $R$. Ceci est absurde, car la fraction $F$ est
bien définie sur les éléments très réguliers donc il existe une
puissance  $\pi^N$ telle que $\pi^N F$ soit un polynôme. On en
déduit que $R^\perp$ est nul et donc que l'homomorphisme de
restriction est surjectif.

\section{Algèbres de Takiff}

\subsection{Rappels sur les algèbres de Takiff}\label{takiff}

Soit $\g$ une algèbre de Lie réductive sur $\C$ de dimension finie.
On suppose que $\g$ est l'algèbre de Lie d'un groupe affine réductif
connexe $G$. Posons $\g_{m}$ pour $m\in \N$ l'algèbre de Lie de
Takiff généralisée \cite{Geo1} définie par
$$\g_{m}=\g\otimes \left(\C[T]/ T^{m+1}\right).$$ C'est une algèbre de Lie algébrique dont un facteur
réductif est isomorphe à $\g$ et le facteur unipotent est
$(\g\otimes T) \oplus \cdots \oplus (\g\otimes T^m)$. Concrètement
on a $$[X\otimes T^{i}, Y\otimes T^{j}]= [X, Y]\otimes T^{i+j}$$si
$i+j \leq  m$ et $0$ sinon. Ces algèbres sont étudiées dans
\cite{Rais-Tauvel} et \cite{Geo1}. On notera $G_m$ un groupe
affine connexe d'algèbre de Lie $G_m$. Ce groupe a un facteur réductif isomorphe à $G$.\\

L'algèbre des invariants $S[\g_m]^{\g_m}$ est une algèbre de
polynômes décrite dans  \cite{Rais-Tauvel} et engendrée par les
dérivées d'ordre inférieur à $m$ des invariants élé\-men\-taires de $S[\g]^\g.$\\

\subsection{L'homomorphisme de restriction pour les algèbres de
Lie algébriques}\label{res}

Pour toute algèbre de Lie $\g$ algébrique, c'est à dire l'algèbre de
Lie d'un groupe algébrique affine connexe $G$ (non nécessairement
réductif) M. Duflo a défini un \ho\ de restriction généralisant
celui de Chevalley. Pour $f \in \g^{*}$ notons $\gf$ le noyau de la
forme bilinéaire alternée associée à $f$, c'est aussi l'annulateur
de $f$ pour l'action coadjointe.  Lorsque $f$ est régulier  $\gf$
est abélien et on note $\sf$ son tore maximal. On dira que $f$ est
générique si $\sf$ est de dimension maximal (parmi les $f$
réguliers). Tous ces tores génériques sont conjugués, ce sont les
sous-algèbres de Cartan-Duflo. Notons $\s$ un tel tore et $\go$ son
 centralisateur dans $\g$. Notons $M'$ le normalisateur de $\s$ dans $G$.
Alors l'\ho\ de restriction généralisant celui de Chevalley dans le
cas réductif est l'application naturelle de $S[\g]^\g$ dans
$S[\go]^{\go}$ déduite de la décomposition $$\g=\go \oplus [\s,
\g].$$C'est une application injective à valeurs dans $S[\go]^{M'}$.
Notons $S(\g)$ le corps des fractions de $S[\g]$ alors on a
l'égalité suivante $S(\g)^{G}\approx S(\go)^{M'}$. Cette égalité
résulte du fait que pour $f\in \go^*$ générique, on a $G\cdot f\cap
\go^*=M'\cdot f$. En général l'\ho\ de restriction n'est pas
surjectif.

\subsection{Description de l'image par la restriction dans le cas Takiff}

Dans cet article on va s'intéresser, lorsque $\g$ est réductive,  à
cet \ho\ de restriction de restriction dans
le cas des algèbres $\g_1=\g\oplus \g\otimes T$, dite algèbre de Takiff.\\

Il existe une forme bilinéaire invariante non dégénérée (mais non
positive) sur $\g_1$ donnée par $$\langle X\otimes T^{\epsilon},
Y\otimes T^{1-\epsilon}\rangle=K(X, Y)$$ avec $K$ la forme de
Killing sur $\g$ (si $\g$ est réductive on prolongera la forme de
Killing en une forme non dégénérée sur $\g$). On identifiera alors
$\g_1$ avec son dual. Pour $f\in \g_1^{*}$ on notera $x_{f}$
l'élément correspondant dans $\g_1$; la sous-algèbre $\gf$ est égale
alors au centralisateur de $x_f$.

Par ailleurs  $x_{f}=x_{o}+x_{1}\otimes T$ est régulier si et
seulement si $x_{o}$ est régulier dans $\g$ \cite{Rais-Tauvel}. Pour
$f$ générique on aura $x_{o}$ générique dans $\g$ c'est à dire
semi-simple régulier.\\

On en déduit que les sous-algèbres de Cartan-Duflo sont conjuguées à
$\h \oplus \h\otimes T$ où $\h$ désigne une sous-algèbre de Cartan
de $\g$. Par conséquent le normalisateur
 dans $G_{1}$ s'identifie au normalisateur dans $G$ c'est à dire le groupe de Weyl $W$.
L'\ho\ de restriction à la Chevalley de $$S[\g_1]^{\g_1} \rightarrow
S[\h_{1}]^{W}$$est donc un isomorphisme au niveau des fractions
rationnelles. L'algèbre  $S[\h_{1}]^{W}$ n'est pas une algèbre de
polynômes, contrairement à $S[\h]^W$. L'homomorphisme de restriction
dans cette situation n'est pas en général surjectif\footnote{Notons
toutefois qu'un résultat de Joseph, assure que l'homomorphisme de
restriction est surjectif de $S[\g_1]^{\g}$ (action diagonale) sur
$S[\h_{1}]^{W}$. Cette restriction n'est cependant pas injective.
}.\\

Les racines de la décomposition de $\g_{1}$ sous l'action de
$\h_{1}$ sont les extensions par $0$ des racines $\A\in \D(\g, \h)$
à $\h_{1}$. On les note encore $\hat{\A}$. On a alors pour tout
$H_{1} +H_{2}\otimes T \in \h_{1}$, $$\langle H_\A\otimes T,\; H_{1}
+H_{2}\otimes T \rangle=\hat{\A}(H_{1}
 +H_{2}\otimes T)=\alpha(H_1),$$avec $H_\A$ la coracine de $\A$ dans $\h$. La coracine
 de $\hat{\A}$ est donc $H_\A\otimes T$ et on a $$\hat{\A}(H_\A\otimes T)=0,$$ c'est à dire que
 dans notre situation les racines sont de carré nul. Dans
\cite{kac}, pour le traitement de l'isomorphisme de Chevalley dans
le cas des algèbres de Kac-Moody ou dans le cas des super-algèbres,
l'existence de racines de carré nul ajoute, dans la description de
l'image des invariants par l'\ho\ de restriction,  une condition sur
les dérivées. Plus précisément on a ici le résultat suivant.\\

Pour $X \in \g_{1}$ notons par $\delta(X)$ la dérivation dans
$S[\g_{1}]$ donnée sur $\g_{1}$ par $$\delta(X)(Y)=\langle X,
Y\rangle.$$ Remarquons alors que  $\delta(H_\A\otimes T)$ correspond
à l'opérateur  $\hat{\A}$ agissant par dérivation.  Pour $\A \in \D$
notons encore par $r_{\A}$ l'action diagonale de la réflexion de $W$
associée à $\A$.

\begin{theo}\label{image}
Soit $p \in S[\h_{1}]$. Alors $p$ est dans l'image de l'\ho\ de
restriction de Chevalley si et seulement si les deux conditions
suivantes sont satisfaites :

1- $p$ est invariant par la $r_{\A}$ pour tout $\A \in \D$

2- $(H_\A\otimes T)^{n}$ divise $\delta^{n}(H_\A\otimes T)p$ pour
tout $n \in \N$.
\end{theo}

\noindent \textit{Preuve:} Le résultat est vrai pour
$\g=\mathrm{sl}(2)$ car on a alors
$$\Res S[\g_{1}]^{\g_{1}}=\C[(H_\A\otimes T)^{2},H_\A (H_\A\otimes T)].$$

On montre d'abord que les conditions sont nécessaires: \\

La première condition est claire d'après le paragraphe précédent. La
deuxième se montre par réduction classique au cas $\mathrm{sl}(2)$.
Soit $h_{o} \in \h$ tel que l'on ait $\A(h_{o})= 0$ et $\beta(
h_{o})\neq 0$ pour $\beta$ racine non proportionnelle à $\A$. Le
centralisateur de $h_{o}$ dans $\g_{1}$ est la somme directe (comme
idéaux) d'une partie abélienne et d'une algèbre de Takiff associée à
un $\mathrm{sl}(2)$ triplet. Comme l'\ho\ de restriction se
factorise par ce
centralisateur, on en déduit la condition~2 du théorème.\\

Montrons que cette condition est suffisante. On raisonne par
récurrence sur la dimension de $\g$.\\

Soit $p$ vérifiant les deux conditions. Alors $p$ s'étend en une
fraction $F$ qui est $G_{1}$-invariante et définie sur les éléments
génériques; c'est toujours la partie facile du théorème de
restriction résultant du fait que les $G_1$-orbites génériques ont
pour traces des $W$-orbites (action diagonale). Sans perte de
généralité on peut supposer que $\g$ est en fait semi-simple.
Montrons que $F$ est partout définie, ce qui montrera que $F$ est un
polynôme. Soit $f\in \g_{1}^{*}$ et notons $x_{f}$ l'élément
correspondant dans $\g_{1}$. Soit $x_{f}=x_{s}+x_{u}$ la
décomposition de Jordan de $x_{f}$. Supposons que l'on ait $x_{s}
\neq 0$. Comme tous les éléments semi-simple de $\g_{1}$ sont
conjugué par $G_{1}$ à un élément de $\h$, on peut supposer sans
perte de généralité que l'on a $x_{s} \in \h$. Notons $\mathfrak{m}$
le centralisateur de
 $x_{s}$ dans $\g$, il vient que le centralisateur dans $\g_1$ vaut
 $\m_{1} =\mathfrak{m}\oplus \mathfrak{m}\otimes T$.
L'algèbre $\mathfrak{m}$ est réductive de dimension strictement
inférieure
 à celle de $\g$ et on a
$x_{f} \in \m_{1}$. Notons $F_{\m_1}$ la restriction de $F$ à
$\m_{1}$, qui est bien définie car on a $\h_{1} \subset \m_{1}$.
L'hypothèse de récurrence assure que  $F_{\m_1}$ est un polynôme,
car les conditions énoncées portent sur les objets liés à $\m$. Or
l'application naturelle de
$$G_{1}\times (\m_{1})^{*} \rightarrow \g_{1}^{*}$$ est une submersion en $f$ donc $F$ est bien
définie en $f$. Reste à étudier le cas des éléments nilpotents. Cet
ensemble vaut $N_{\g} + \g\otimes T$ où $N_{\g}$ désigne les
éléments nilpotents de $\g$. Cet ensemble ne constitue pas une
sous-variété de codimension $1$ sauf dans le cas $\mathrm{sl}(2)$,
que l'on a déjà étudié. D'après le principe de Hartogs $F$ est donc
un polynôme. $\Box$
\\

\begin{cor}
Comme dans le théorème de Harish-Chandra, on déduit de la
proposition un résultat analogue pour le centre de
 l'algèbre enveloppante de~$\g_{1}$.
\end{cor}

\begin{rem}
Pour l'algèbre de Takiff généralisée $\mathrm{sl}(2)_2$ il n'est pas
difficile de voir l'image par l'homomorphisme de restriction est
engendrée par $$(H\otimes T^2)^2,\quad   (H\otimes T^2)(H\otimes T),
\quad (H\otimes T)^2 +2 H(H\otimes T^2).$$Les conditions énoncées
dans le théorème \ref{image} ne sont donc pas suffisantes dans ce
cas pour décrire l'image.

\end{rem}

\subsection{Cas des paires symétriques de
Takiff}

Soit $(\g, \theta)$ une paire sym\'etrique r\'eductive et d\'efinie
sur $\mathbb{C}$. L'involution $\theta$ est obtenue par extension
\`a $\mathbb{C}$ d'une involution de Cartan sur une forme r\'eelle
de $\g$. On notera $\g=\k\oplus \p$ la d\'ecomposition relativement
\`a $\theta$. On appelle $\theta$-rang la dimension d'un sous-espace de Cartan
dans $\p$.\\

On \'etend l'involution $\theta$ \`a $\g_1=\g\oplus \g\otimes T$ de
mani\`ere \'evidente; on note $\theta_1$ une telle extension et
\'ecrit  $\g_1=\k_1\oplus \p_1$ pour la d\'ecomposition associ\'ee.
Alors la paire $(\g_1, \sigma_1)$ sera appel\'ee paire sym\'etrique
de Takiff. Plus g\'en\'eralement on pourra \'etendre $\theta$ \`a
$$\g_n=\g\oplus \cdots \oplus \g\otimes T^n=\k_n\oplus \p_n,$$ on
parlera alors de paires sym\'etriques de Takiff
g\'en\'eralis\'ees.\\

L'involution $\theta$ se remonte \`a $G$ (groupe r\'eductif affine
connexe). On note $K_\theta$ le sous-groupe des  points fixes et $K$
sa composante irr\'eductible.  Les groupes $K$ et $K_\theta$
diffèrent d'un groupe fini d'éléments d'ordre $2$
(\cite{Kost-Rallis} $\S$ I.1 prop. 1).\\

Les groupes $K$ et $K_\theta$ agissent sur $S[\p]$, l'alg\`ebre
sym\'etrique de $\p$. D'apr\`es \cite{Kost-Rallis} $\S$ I.4 prop. 10
on a
$$S[\p]^\k=S[\p]^K=S[\p]^{K_\theta}.$$

On généralise sans peine les résultats de Rais-Tauvel sur la
description des invariants $S[\p_1]^{k_1}$ (resp. $S[\p_n]^{\k_n}$),
grâce aux résultats de
Kostant-Rallis sur les espaces symétriques~\cite{Kost-Rallis}.\\

La description des invariants par le groupe $K$ ($K_\theta$) et
l'existence d'une transversale de Kostant pour les
$K_\theta$-orbites (\cite{Kost-Rallis} \S II.3), assure l'extension
des arguments et les résultats de \cite{Rais-Tauvel}, notamment le
lemme 3.5 et le théorème 4.5. En particulier l'alg\`ebre des
invariants $S[\p_n]^{k_n}$ forme une algèbre de polynômes engendrée
par les dérivées des générateurs de l'algèbre $S[\p]^{k}$.

Remarquons que les  paires de Takiff de $\theta$-rang $1$, ont déjà
été rencontr\'ees dans \cite{SMF-rad} \S 2.4.2. \\

Dans le cas des paires symétriques de Takiff sur $\mathbb{C}$, les
arguments dévelop\-pés dans cet article assurent sans peine
l'extension du théorème \ref{image} au cas des espaces symétriques
de Takiff. Principalement on remplace les sous-alg\`ebres de Cartan
par le sous-espace de Cartan et $\mathrm{sl}(2)$ est remplac\'e par
les paires sym\'etriques de $\theta$-rang $1$.\\

 En effet la r\'eduction
proc\`ede du m\^eme principe, on raisonne par r\'ecurrence sur le
rang de la paire sym\'etrique de Takiff. Comme dans la preuve du
théorème \ref{image} on consid\`ere le centralisateur d'un
\'el\'ement $H_o $ d'un sous-espace de Cartan v\'erifiant
$\alpha(H_o)= 0$ (pour $\alpha$ une racine r\'eduite) et
$\beta(H_o)\neq 0$ pour $\beta$ toute racine non proportionnelle \`a
$\alpha$. Le centralisateur de $H_o$ nous ram\`ene \`a une paire
sym\'etrique de Takiff de rang $1$. L'image par l'homomorphisme de
restriction dans ce cas là a encore pour image (avec les notations
du théorème~\ref{image})
$$(H_\A\otimes T)^{2}, \quad H_\A (H_\A\otimes T).$$

La r\'ecurrence peut donc s'amorcer et le reste des arguments est
identique, y compris l'argument sur la co-dimension des éléments
nilpotents car  pour une paire sym\'etrique réductive sur
$\mathbb{C}$ les \'el\'ements nilpotents forment encore une
vari\'et\'e de codimension égale au  rang (\cite{Kost-Rallis} \S I.4
th. 3 ).

\end{document}